\documentclass[a4paper]{article}
\usepackage[T1]{fontenc}
\usepackage[latin1]{inputenc}
\usepackage{lmodern}
\usepackage[french]{babel}
\usepackage{epsfig}
\usepackage{colortbl}
\usepackage{hyperref}
\usepackage{graphicx}
\usepackage{lscape}
\usepackage{amsfonts}
\usepackage{amssymb}
\usepackage{amsmath}
\title{Nombre de r\'{e}sidus quadratiques \\d'un nombre entier \\inf\'{e}rieurs \`{a} sa moiti\'{e}}
\author{Denise Vella-Chemla}
\date{26.8.2016}

\begin{document}
\maketitle

\noindent On souhaiterait d\'{e}montrer qu'on peut \'{e}tablir le caract\`{e}re de primalit\'{e} d'un entier $n$ impair en comptant le nombre (qu'on note $R_b(n)$) de ses r\'{e}sidus quadratiques (non nuls, on le sp\'{e}cifie une fois pour toutes) qui sont inf\'{e}rieurs ou \'{e}gaux \`{a} $n/2$.

\noindent Plus pr\'{e}cis\'{e}ment, on induit de comptages effectu\'{e}s pour les nombres impairs jusqu'\`{a} 100 l'hypoth\`{e}se suivante :

\noindent \begin{textit}{(H) Si le nombre $R_b(n)$ des r\'{e}sidus quadratiques d'un entier $n$ inf\'{e}rieurs \`{a} $n/2$ est sup\'{e}rieur \`{a} $n/4$, $n$ est premier ; sinon, $n$ est compos\'{e}.}\end{textit}

\noindent Cette hypoth\`{e}se s'\'{e}crit :\\
$$
\begin{array}{l}
(H) \;\;\;\forall p, \;p \;impair \;et \;p \;\leq \;3,\\
\hspace{2cm}R_b(p) \;= \;\# \;\left\{y \;tels \;que \;\exists x \in \mathbb{N}^{\times}, \;\exists k \in \mathbb{N}, \;x^2-kp-y=0 \;avec \;0 \;< \;y \;\leq \;\displaystyle\frac{p-1}{2} \;\right\} > \;\displaystyle\frac{n}{4}\\
\hspace{1cm}\iff\\
\hspace{2cm}p \;premier
\end{array}
$$

\noindent On teste notre hypoth\`{e}se en la programmant : on ne peut la tester tr\`{e}s loin car l'\'{e}l\'{e}vation au carr\'{e} d\'{e}passe vite les limites des entiers (la limite des ``unsigned long int'' en C++ ne permet de tester l'hypoth\`{e}se que pour des entiers inf\'{e}rieurs \`{a} $300\;000$ car on utilise un programme simple). On r\'{e}alise en programmant que les nombres premiers $p$ de la forme $4k+1$ ont leur nombre de ``petits r\'{e}sidus quadratiques'' qui est \'{e}gal \`{a} $\lfloor p/4 \rfloor$ tandis que les nombres premiers $p$ de la forme $4k+3$ ont le nombre en question strictement sup\'{e}rieur \`{a} $\lfloor p/4 \rfloor$. 

\noindent On va utiliser la notation $x \;R \;m$ pour exprimer que $x$ est un r\'{e}sidu quadratique de $m$ et la notation $x \;N \;m$ pour exprimer que $x$ n'est pas un r\'{e}sidu quadratique de $m$.

\noindent On rappelle la d\'{e}finition de \begin{textit}{x est un r\'{e}sidu quadratique de m}\end{textit} (not\'{e}e $x \;R \;m$) :
$$x \;R \;m \;\iff \;\exists y \;tel \;que \;x^2 \;\equiv \;y \;(mod \;m).$$

\noindent On a :
\begin{itemize}
\item{$\forall x, y, m, \;x \;R \;m \;et \;y \;R \;m \;\implies \;xy \;R \;m$}

\item{$\forall x, y, m, \;x \;N \;m \;et \;y \;N \;m \;\implies \;xy \;R \;m$}

\item{$\forall x, y, m, \;x \;R \;m \;et \;y \;N \;m \;\implies \;xy \;N \;m$}
\end{itemize}

\vspace{0.5cm}
\noindent En annexe, sont fournis les r\'{e}sidus quadratiques des nombres entiers impairs de 3 \`{a} 51 inf\'{e}rieurs \`{a} leur moiti\'{e} ainsi que leur nombre.

\noindent On rappelle qu'on peut d\'{e}duire le caract\`{e}re de r\'{e}siduosit\'{e} quadratique d'un nombre $a$ premier \`{a} $p$ un nombre premier impair par un simple calcul de puissance (cf. paragraphe 106 des Recherches arithm\'{e}tiques de Gauss) :
$$
\left(\frac{a}{p}\right)=a^{\frac{p-1}{2}} \;(mod \;p)
$$

  \noindent On rappelle que $-1$ est r\'{e}sidu quadratique de tout nombre premier de la forme $4k+1$ et n'est pas r\'{e}sidu quadratique de tout nombre premier de la forme $4k+3$ (cf. paragraphe 109 des Recherches arithm\'{e}tiques de Gauss).

\noindent Ces particularit\'{e}s du caract\`{e}re de r\'{e}siduosit\'{e} quadratique de $-1$ aux nombres premiers impairs ont pour cons\'{e}quence que si $p$ est un nombre premier impair de la forme $4k+1$ alors $r \;R\; p \iff -r \;R\; p$ tandis que si $p$ est un nombre premier impair de la forme $4k+3$ alors $r \;R\; p \iff -r \;N\; p$.

\noindent Pour d\'{e}montrer notre hypoth\`{e}se, il faudrait d\'{e}montrer :

 1) qu'elle est vraie pour les nombres premiers de la forme $4k+1$ ;

 2) qu'elle est vraie pour les nombres premiers de la forme $4k+3$ ;

 3) qu'elle est vraie par \'{e}l\'{e}vation d'un nombre premier $p$ \`{a} la puissance $k$ ;

 4) qu'elle est vraie par multiplication de deux nombres premiers simples ;

 5) qu'elle est vraie par multiplication de deux puissances de nombres premiers (incluant peut-\^{e}tre (4)).

\vspace{0.3cm}

\noindent 1) Pour les nombres premiers $p$ de la forme $4k+1$, le nombre de r\'{e}sidus quadratiques de $p$ inf\'{e}rieurs ou \'{e}gaux \`{a} $p/2$ est trivialement \'{e}gal \`{a} $\displaystyle\frac{p-1}{4}$ car $r$ et $-r$ sont syst\'{e}matiquement soit tous deux r\'{e}sidus quadratiques de $p$ soit tous deux non-r\'{e}sidus quadratiques de $p$ et parce que seul l'un des deux dans chaque couple est inf\'{e}rieur ou \'{e}gal \`{a} $p/2$.

\noindent 2) Pour les nombres premiers $p$ de la forme $4k+3$, Dirichlet a d\'{e}montr\'{e} qu'il y a davantage de petits r\'{e}sidus quadratiques de $p$ que de petits non-r\'{e}sidus quadratiques de $p$ en utilisant l'analyse infinit\'{e}simale\footnote{Je ne parviens pas \`{a} trouver ce r\'{e}sultat de Dirichlet dont les r\'{e}f\'{e}rences possibles sont \textit{Applications de l'analyse infinit\'{e}simale \`{a} la th\'{e}orie des nombres}, Journal de Crelle, 1839, vol. 19, p.324-369 ou vol. 21, p.1-12 ou p.134-155.} :
$$
\sum\limits_{n=1}^{n=\frac{p-1}{2}}\left(\displaystyle\frac{n}{p}\right) > 0.
$$

\noindent Dans la suite, l'adjectif \textit{petit} signifie \textit{inf\'{e}rieur ou \'{e}gal \`{a} $p/2$} et \textit{grand} signifie \textit{sup\'{e}rieur strictement \`{a} $p/2$}.

\noindent Dans [1], on trouve les relations invariantes suivantes ; notons :

\begin{itemize}
\item{$\sum R$ la somme des r\'{e}sidus quadratiques de $p$ un nombre premier impair,}
\item{$\sum N$ la somme de ses non-r\'{e}sidus quadratiques,}
\item{$\sum R_b$ la somme de ses ``petits'' r\'{e}sidus quadratiques,}
\item{$\sum N_b$ la somme de ses ``petits'' non-r\'{e}sidus quadratiques,}
\item{$\sum R_h$ la somme de ses ``grands'' r\'{e}sidus quadratiques,}
\item{$\sum N_h$ la somme de ses ``grands'' non-r\'{e}sidus quadratiques,}
\item{$R_b$ le nombre de ses petits r\'{e}sidus quadratiques ;}
\item{$N_b$ le nombre de ses petits non-r\'{e}sidus quadratiques.}
\end{itemize}

\vspace{0.5cm}
\noindent On a, si $p \equiv 7 \;(mod \;8)$ :\\
$$
\left\{
\begin{array}{l}
\sum R_b = \sum N_b.\\\\
\displaystyle\frac{\sum N-\sum R}{p} = R_b-N_b.
\end{array}
\right.
$$

\noindent On a, si $p \equiv 3 \;(mod \;8)$ :\\
$$
\left\{
\begin{array}{l}
\sum R - \sum N = \sum R_b- \sum N_b.\\\\
\displaystyle\frac{3(\sum N-\sum R)}{p} = R_b-N_b.
\end{array}
\right.
$$
\newpage
\noindent Concernant les nombres premiers de la forme $p=8k+7$, on r\'{e}alise par le calcul, m\^{e}me si on ne sait pas le d\'{e}montrer, que la somme des petits r\'{e}sidus quadratiques, qui est \'{e}gale selon l'un des r\'{e}sultats de Victor-Am\'{e}d\'{e}e Lebesgue \`{a} la somme des petits non-r\'{e}sidus quadratiques, est \'{e}gale \`{a} :$$\sum R_b=\displaystyle\frac{(p-1)(p+1)}{16}.$$

\noindent Il y a s\^{u}rement de tr\`{e}s nombreuses fonctions $f$ qui sont telles que $f(7)=3,f(23)=33,f(31)=60,f(47)=138,f(71)=315,f(79)=390,f(103)=663$ et $f(9967)=6208818$ mais $f(p)=\displaystyle\frac{(p-1)(p+1)}{16}$ semble pertinente dans ce contexte.

\noindent On a ainsi deux ensembles, l'ensemble des petits r\'{e}sidus quadratiques et l'ensemble des petits non-r\'{e}sidus quadratiques, dont les sommes des \'{e}l\'{e}ments sont \'{e}gales \`{a} $\displaystyle\frac{(p-1)(p+1)}{16}$, et dont on sait qu'ils sont de cardinaux diff\'{e}rents. On sait \'{e}galement que $4$ est toujours \'{e}l\'{e}ment de l'ensemble des petits r\'{e}sidus quadratiques. Ces diff\'{e}rents \'{e}l\'{e}ments ont peut-\^{e}tre pour cons\'{e}quence qu'il y a forc\'{e}ment plus de petits r\'{e}sidus quadratiques que de petits non-r\'{e}sidus quadratiques.

\noindent Pour illustrer l'id\'{e}e de la sup\'{e}riorit\'{e} en nombre des petits r\'{e}sidus quadratiques, on peut pr\'{e}senter la multiplication modulaire modulo 7 sur 2 tables, l'une dans laquelle les nombres sont comme habituellement dans l'ordre croissant, l'autre dans laquelle les r\'{e}sidus quadratiques sont \'{e}num\'{e}r\'{e}s avant les non-r\'{e}sidus quadratiques.

$$
\begin{tabular}{c|ccc|ccc|}
  &1&2&3&4&5&6\\
  \hline
  1&\cellcolor{cyan}1&\cellcolor{cyan}2&3&\cellcolor{cyan}4&5&6\\
  2&\cellcolor{cyan}2&\cellcolor{cyan}4&6&\cellcolor{cyan}1&3&5\\
  3&3&6&\cellcolor{cyan}2&5&\cellcolor{cyan}1&\cellcolor{cyan}4\\
  \hline
  4&\cellcolor{cyan}4&\cellcolor{cyan}1&5&\cellcolor{cyan}2&6&3\\
  5&5&3&\cellcolor{cyan}1&6&\cellcolor{cyan}4&\cellcolor{cyan}2\\
  6&6&5&\cellcolor{cyan}4&3&\cellcolor{cyan}2&\cellcolor{cyan}1\\
  \hline
  \end{tabular}
$$\\

$$
\begin{tabular}{c|ccc|ccc|}
  &1&2&4&3&5&6\\
  \hline
  1&\cellcolor{cyan}1&\cellcolor{cyan}2&\cellcolor{cyan}4&3&5&6\\
  2&\cellcolor{cyan}2&\cellcolor{cyan}4&\cellcolor{cyan}1&6&3&5\\
  4&\cellcolor{cyan}4&\cellcolor{cyan}1&\cellcolor{cyan}2&5&6&3\\
  \hline
  3&3&6&5&\cellcolor{cyan}2&\cellcolor{cyan}1&\cellcolor{cyan}4\\
  5&5&3&6&\cellcolor{cyan}1&\cellcolor{cyan}4&\cellcolor{cyan}2\\
  6&6&5&3&\cellcolor{cyan}4&\cellcolor{cyan}2&\cellcolor{cyan}1\\
  \hline
  \end{tabular}
$$\\

\noindent Sur la table de multiplication modulaire modulo 23, la sup\'{e}riorit\'{e} en nombre des petits r\'{e}sidus quadratiques sur les petits non-r\'{e}sidus quadratiques apparaît clairement dans le quart haut-gauche de la table. Ce sont peut-\^{e}tre les permutations qui perturbent l'ordre habituel sur les entiers lorsqu'on les consid\`{e}re comme \'{e}tant des carr\'{e}s modulaires qui a pour cons\'{e}quence la loi qui nous semble \^{e}tre toujours v\'{e}rifi\'{e}e (ou du moins jusqu'\`{a} $300\;000$).

\begin{scriptsize}
$$
  \begin{tabular}{c|ccccccccccc|ccccccccccc|}
    &1 & 2 & 3 & 4 & 6 & 8 & 9 & 12 & 13 & 16 & 18 & 5 & 7 & 10 & 11 & 14 & 15 & 17 & 19 & 20 & 21 & 22 \\
    \hline
1&\cellcolor{cyan}1 & \cellcolor{cyan}2 & \cellcolor{cyan}3 & \cellcolor{cyan}4 & \cellcolor{cyan}6 & \cellcolor{cyan}8 & \cellcolor{cyan}9 & 12 & 13 & 16 & 18 & \cellcolor{cyan}5 & \cellcolor{cyan}7 & \cellcolor{cyan}10 & \cellcolor{cyan}11 & 14 & 15 & 17 & 19 & 20 & 21 & 22 \\
2&\cellcolor{cyan}2 & \cellcolor{cyan}4 & \cellcolor{cyan}6 & \cellcolor{cyan}8 & 12 & 16 & 18 & \cellcolor{cyan}1 & \cellcolor{cyan}3 & \cellcolor{cyan}9 & 13 & \cellcolor{cyan}10 & 14 & 20 & 22 & \cellcolor{cyan}5 & \cellcolor{cyan}7 & \cellcolor{cyan}11 & 15 & 17 & 19 & 21 \\
3&\cellcolor{cyan}3 & \cellcolor{cyan}6 & \cellcolor{cyan}9 & 12 & 18 & \cellcolor{cyan}1 & \cellcolor{cyan}4 & 13 & 16 & \cellcolor{cyan}2 & \cellcolor{cyan}8 & 15 & 21 & \cellcolor{cyan}7 & \cellcolor{cyan}10 & 19 & 22 & \cellcolor{cyan}5 & \cellcolor{cyan}11 & 14 & 17 & 20 \\
4&\cellcolor{cyan}4 & \cellcolor{cyan}8 & 12 & 16 & \cellcolor{cyan}1 & \cellcolor{cyan}9 & 13 & \cellcolor{cyan}2 & \cellcolor{cyan}6 & 18 & \cellcolor{cyan}3 & 20 & \cellcolor{cyan}5 & 17 & 21 & \cellcolor{cyan}10 & 14 & 22 & \cellcolor{cyan}7 & \cellcolor{cyan}11 & 15 & 19 \\
6&\cellcolor{cyan}6 & 12 & 18 & \cellcolor{cyan}1 & 13 & \cellcolor{cyan}2 & \cellcolor{cyan}8 & \cellcolor{cyan}3 & \cellcolor{cyan}9 & \cellcolor{cyan}4 & 16 & \cellcolor{cyan}7 & 19 & 14 & 20 & 15 & 21 & \cellcolor{cyan}10 & 22 & \cellcolor{cyan}5 & \cellcolor{cyan}11 & 17 \\
8&\cellcolor{cyan}8 & 16 & \cellcolor{cyan}1 & \cellcolor{cyan}9 & \cellcolor{cyan}2 & 18 & \cellcolor{cyan}3 & \cellcolor{cyan}4 & 12 & 13 & \cellcolor{cyan}6 & 17 & \cellcolor{cyan}10 & \cellcolor{cyan}11 & 19 & 20 & \cellcolor{cyan}5 & 21 & 14 & 22 & \cellcolor{cyan}7 & 15 \\
9&\cellcolor{cyan}9 & 18 & \cellcolor{cyan}4 & 13 & \cellcolor{cyan}8 & \cellcolor{cyan}3 & 12 & 16 & \cellcolor{cyan}2 & \cellcolor{cyan}6 & \cellcolor{cyan}1 & 22 & 17 & 21 & \cellcolor{cyan}7 & \cellcolor{cyan}11 & 20 & 15 & \cellcolor{cyan}10 & 19 & \cellcolor{cyan}5 & 14 \\
12&12 & \cellcolor{cyan}1 & 13 & \cellcolor{cyan}2 & \cellcolor{cyan}3 & \cellcolor{cyan}4 & 16 & \cellcolor{cyan}6 & 18 & \cellcolor{cyan}8 & \cellcolor{cyan}9 & 14 & 15 & \cellcolor{cyan}5 & 17 & \cellcolor{cyan}7 & 19 & 20 & 21 & \cellcolor{cyan}10 & 22 & \cellcolor{cyan}11 \\ 
13&13 & \cellcolor{cyan}3 & 16 & \cellcolor{cyan}6 & \cellcolor{cyan}9 & 12 & \cellcolor{cyan}2 & 18 & \cellcolor{cyan}8 & \cellcolor{cyan}1 & \cellcolor{cyan}4 & 19 & 22 & 15 & \cellcolor{cyan}5 & 21 & \cellcolor{cyan}11 & 14 & 17 & \cellcolor{cyan}7 & 20 & \cellcolor{cyan}10 \\
16&16 & \cellcolor{cyan}9 & \cellcolor{cyan}2 & 18 & \cellcolor{cyan}4 & 13 & \cellcolor{cyan}6 & \cellcolor{cyan}8 & \cellcolor{cyan}1 & \cellcolor{cyan}3 & 12 & \cellcolor{cyan}11 & 20 & 22 & 15 & 17 & \cellcolor{cyan}10 & 19 & \cellcolor{cyan}5 & 21 & 14 & \cellcolor{cyan}7 \\
18&18 & 13 & \cellcolor{cyan}8 & \cellcolor{cyan}3 & 16 & \cellcolor{cyan}6 & \cellcolor{cyan}1 & \cellcolor{cyan}9 & \cellcolor{cyan}4 & 12 & \cellcolor{cyan}2 & 21 & \cellcolor{cyan}11 & 19 & 14 & 22 & 17 & \cellcolor{cyan}7 & 20 & 15 & \cellcolor{cyan}10 & \cellcolor{cyan}5 \\
\hline
5&\cellcolor{cyan}5 & \cellcolor{cyan}10 & 15 & 20 & \cellcolor{cyan}7 & 17 & 22 & 14 & 19 & \cellcolor{cyan}11 & 21 & \cellcolor{cyan}2 & 12 & \cellcolor{cyan}4 & \cellcolor{cyan}9 & \cellcolor{cyan}1 & \cellcolor{cyan}6 & 16 & \cellcolor{cyan}3 & \cellcolor{cyan}8 & 13 & 18 \\
7&\cellcolor{cyan}7 & 14 & 21 & \cellcolor{cyan}5 & 19 & \cellcolor{cyan}10 & 17 & 15 & 22 & 20 & \cellcolor{cyan}11 & 12 & \cellcolor{cyan}3 & \cellcolor{cyan}1 & \cellcolor{cyan}8 & \cellcolor{cyan}6 & 13 & \cellcolor{cyan}4 & 18 & \cellcolor{cyan}2 & \cellcolor{cyan}9 & 16 \\
10&\cellcolor{cyan}10 & 20 & \cellcolor{cyan}7 & 17 & 14 & \cellcolor{cyan}11 & 21 & \cellcolor{cyan}5 & 15 & 22 & 19 & \cellcolor{cyan}4 & \cellcolor{cyan}1 & \cellcolor{cyan}8 & 18 & \cellcolor{cyan}2 & 12 & \cellcolor{cyan}9 & \cellcolor{cyan}6 & 16 & \cellcolor{cyan}3 & 13 \\ 
11&\cellcolor{cyan}11 & 22 & \cellcolor{cyan}10 & 21 & 20 & 19 & \cellcolor{cyan}7 & 17 & \cellcolor{cyan}5 & 15 & 14 & \cellcolor{cyan}9 & \cellcolor{cyan}8 & 18 & \cellcolor{cyan}6 & 16 & \cellcolor{cyan}4 & \cellcolor{cyan}3 & \cellcolor{cyan}2 & 13 & \cellcolor{cyan}1 & 12 \\
14&14 & \cellcolor{cyan}5 & 19 & \cellcolor{cyan}10 & 15 & 20 & \cellcolor{cyan}11 & \cellcolor{cyan}7 & 21 & 17 & 22 & \cellcolor{cyan}1 & \cellcolor{cyan}6 & \cellcolor{cyan}2 & 16 & 12 & \cellcolor{cyan}3 & \cellcolor{cyan}8 & 13 & \cellcolor{cyan}4 & 18 & \cellcolor{cyan}9 \\
15&15 & \cellcolor{cyan}7 & 22 & 14 & 21 & \cellcolor{cyan}5 & 20 & 19 & \cellcolor{cyan}11 & \cellcolor{cyan}10 & 17 & \cellcolor{cyan}6 & 13 & 12 & \cellcolor{cyan}4 & \cellcolor{cyan}3 & 18 & \cellcolor{cyan}2 & \cellcolor{cyan}9 & \cellcolor{cyan}1 & 16 & \cellcolor{cyan}8 \\
17&17 & \cellcolor{cyan}11 & \cellcolor{cyan}5 & 22 & \cellcolor{cyan}10 & 21 & 15 & 20 & 14 & 19 & \cellcolor{cyan}7 & 16 & \cellcolor{cyan}4 & \cellcolor{cyan}9 & \cellcolor{cyan}3 & \cellcolor{cyan}8 & \cellcolor{cyan}2 & 13 & \cellcolor{cyan}1 & 18 & 12 & \cellcolor{cyan}6 \\
19&19 & 15 & \cellcolor{cyan}11 & \cellcolor{cyan}7 & 22 & 14 & \cellcolor{cyan}10 & 21 & 17 & \cellcolor{cyan}5 & 20 & \cellcolor{cyan}3 & 18 & \cellcolor{cyan}6 & \cellcolor{cyan}2 & 13 & \cellcolor{cyan}9 & \cellcolor{cyan}1 & 16 & 12 & \cellcolor{cyan}8 & \cellcolor{cyan}4 \\
20&20 & 17 & 14 & \cellcolor{cyan}11 & \cellcolor{cyan}5 & 22 & 19 & \cellcolor{cyan}10 & \cellcolor{cyan}7 & 21 & 15 & \cellcolor{cyan}8 & \cellcolor{cyan}2 & 16 & 13 & \cellcolor{cyan}4 & \cellcolor{cyan}1 & 18 & 12 & \cellcolor{cyan}9 & \cellcolor{cyan}6 & \cellcolor{cyan}3 \\
21&21 & 19 & 17 & 15 & \cellcolor{cyan}11 & \cellcolor{cyan}7 & \cellcolor{cyan}5 & 22 & 20 & 14 & \cellcolor{cyan}10 & 13 & \cellcolor{cyan}9 & \cellcolor{cyan}3 & \cellcolor{cyan}1 & 18 & 16 & 12 & \cellcolor{cyan}8 & \cellcolor{cyan}6 & \cellcolor{cyan}4 & \cellcolor{cyan}2 \\
22&22 & 21 & 20 & 19 & 17 & 15 & 14 & \cellcolor{cyan}11 & \cellcolor{cyan}10 & \cellcolor{cyan}7 & \cellcolor{cyan}5 & 18 & 16 & 13 & 12 & \cellcolor{cyan}9 & \cellcolor{cyan}8 & \cellcolor{cyan}6 & \cellcolor{cyan}4 & \cellcolor{cyan}3 & \cellcolor{cyan}2 & \cellcolor{cyan}1 \\
\hline
\end{tabular}
$$
\end{scriptsize}

\noindent Noter les sym\'{e}tries-miroir horizontale et verticale.

\noindent Int\'{e}ressons-nous au cas qu'on avait num\'{e}rot\'{e} (3) ci-dessus dans le cas o\`{u} on parviendrait \`{a} d\'{e}montrer formellement notre hypoth\`{e}se, c'est \`{a} dire le fait que pour les nombres $n$ qui sont des puissances de nombres premiers ($n=p^m$ avec $p$ premier), il y a moins de $n/4$ petits r\'{e}sidus quadratiques qui sont inf\'{e}rieurs \`{a} $n/2$.  

\noindent Pour le calcul du nombre de petits r\'{e}sidus quadratiques de $p^k$ une puissance d'un nombre premier $p$ (i.e. les r\'{e}sidus inf\'{e}rieurs ou \'{e}gaux \`{a} $p/2$), on induit des formules en analysant les premi\`{e}res valeurs calcul\'{e}es par ordinateur pour 3 d'abord, car il semble se comporter diff\'{e}remment des autres, puis pour les puissances des nombres premiers de la forme $4k+1$ et enfin, pour les puissances des nombres premiers de la forme $4k+3$.

\noindent Des valeurs
$$
\begin{array}{ll}
  R_b(9)&=2,\\
  R_b(27)&=6,\\
  R_b(81)&=16,\\
  R_b(243)&=47,\\
  R_b(729)&=138,\\
  R_b(2187)&=412,\\
  R_b(6561)&=1232,\\
  R_b(19683)&=3693,
\end{array}
$$
\noindent on induit la formule suivante pour le nombre de petits r\'{e}sidus quadratiques des puissances de $3$.
$$
R_b(3^k)=\left\lceil\displaystyle\frac{3^{k-1}}{2}\right\rceil+R_b(3^{k-2}).
$$

\noindent Des valeurs
$$
\begin{array}{ll}
  R_b(25)&=5,\\
  R_b(125)&=26,\\
  R_b(625)&=130,\\
  R_b(3125)&=651,\\
  R_b(15625)&=3255,
\end{array}
$$
\noindent pour les puissances de $5$, puis
$$
\begin{array}{ll}
  R_b(169)&=39,\\
  R_b(2197)&=510,\\
  R_b(28561)&=6630,
  \end{array}
$$
\noindent pour les puissances de $13$, on trouve la formule suivante pour le nombre de petits r\'{e}sidus quadratiques des puissances de $p$ pour $p$ de la forme $4k+1$.
$$
R_b(p^k)=kp^{k-1}+R_b(p^{k-2}).
$$

\noindent Ce qu'il est amusant de constater, c'est que cette formule pour les $4k+1$ s'applique \'{e}galement au nombre premier $3$ si ce n'est qu'il faut consid\'{e}rer celui-ci non pas comme un $4k+3$, ce qu'on a coutume de faire habituellement, mais plut\^{o}t comme un $4k+1$ avec $k$ qui vaudrait $\displaystyle\frac{1}{2}$.

\noindent Des valeurs
$$
\begin{array}{ll}
R_b(7)&=2,\\
R_b(49)&=11,\\
R_b(343)&=76,\\
R_b(2401)&=526,\\
R_b(16807)&=3678,
\end{array}
$$
\noindent pour les puissances de $7$, puis
$$
\begin{array}{ll}
R_b(11)&=4,\\
R_b(121)&=29,\\
R_b(1331)&=308,\\
R_b(14641)&=3358,
\end{array}
$$
\noindent pour les puissances de $11$, on induit que le nombre de petits r\'{e}sidus quadratiques des puissances de $p$ pour $p$ de la forme $4k+3$ (sauf 3) est d'un ordre proche de $p$ fois le nombre de petits r\'{e}sidus quadratiques de la puissance de $p$ juste inf\'{e}rieure :
$$
R_b(p^k) \approx pR_b(p^{k-1}).
$$

\noindent Il semble plus judicieux, par un principe de sym\'{e}trie, de s'int\'{e}resser au nombre de grands non-r\'{e}sidus (i.e. sup\'{e}rieurs strictement \`{a} $\displaystyle\frac{n-1}{2}$), qu'on note $N_h(n)$. En effet, les valeurs des $N_h(7^k)$ ou des $N_h(11^k)$ sont :
$$
\begin{array}{ll}
  N_h(7)&=2,\\
  N_h(49=7^2)&=14,\\
  N_h(343=7^3)&=97,\\
  N_h(2401=7^4)&=676,\\\\

  N_h(11)&=4,\\
  N_h(121=11^2)&=34,\\
  N_h(1331=11^3)&=363.
  \end{array}
$$

\noindent L'approximation $N_h(p^k) \approx pN_h(p^{k-1})$ semble de meilleure qualit\'{e} que l'approximation pr\'{e}c\'{e}dente.
  
\noindent Concernant les produits de puissances (point (5) \'{e}voqu\'{e} plus haut), des valeurs suivantes de $N_h(n)$,
$$
\begin{array}{ll}
  N_h(3.7)&=7,\\
  N_h(3^2.7)&=25,\\
  N_h(3.7^2)&=52,\\
  N_h(3^2.7^2&=178,\\\\
  
  N_h(5.7)&=13,\\
  N_h(5^2.7)&=68,\\
  N_h(5.7^2)&=91,\\
  N_h(5^2.7^2&=494,\\\\
  
  N_h(13.17)&=79,\\
  N_h(13^2.17)&=1081,\\
  N_h(13.17^2)&=1399.\\
  \end{array}
$$
\noindent on induit un ordre d'approximation pour $N_h(n)$ ($p<q$) :
$$
N_h(p^m.q^n) \approx p.N_h(p^{m-1}.q^n)
$$

\noindent Concernant les nombres de petits r\'{e}sidus quadratiques dont les valeurs suivent,
$$
\begin{array}{ll}
  R_b(3.7)&=4,\\
  R_b(3^2.7)&=9,\\
  R_b(3.7^2)&=22,\\
  R_b(3^2.7^2&=45,\\\\
  
  R_b(5.7)&=7,\\
  R_b(5^2.7)&=24,\\
  R_b(5.7^2)&=34,\\
  R_b(5^2.7^2&=123,\\\\
  
  R_b(13.17)&=31,\\
  R_b(13^2.17)&=355,\\
  R_b(13.17^2)&=479.\\
  \end{array}
$$
\noindent m\^{e}me si on ne peut trouver de r\`{e}gles quant \`{a} leur progression, on constate que les valeurs \'{e}tudi\'{e}es v\'{e}rifient toujours :
$$
R_b(p^m.q^n) < p.R_b(p^{m-1}.q^n).
$$

\noindent Le nombre de petits r\'{e}sidus quadratiques pour les puissances de nombres premiers de la forme $p=4k+3$ (sauf dans le cas o\`{u} $p=3$) est toujours strictement inf\'{e}rieur \`{a} $\displaystyle\frac{p^k-1}{4}$ car tous les multiples de $p$ ne peuvent \^{e}tre r\'{e}sidus quadratiques des puissances successives de $p$, ce qui r\'{e}duit d'autant le nombre de petits r\'{e}sidus quadratiques.

\noindent On peut maintenant d\'{e}monter le m\'{e}canisme \`{a} l'oeuvre pour les produits, ce m\'{e}canisme ayant pour cons\'{e}quence une tr\`{e}s grande redondance des carr\'{e}s obtenus, qui r\'{e}duit d'autant leur nombre, le rendant toujours inf\'{e}rieur au quart du nombre $n$ consid\'{e}r\'{e}.

\noindent Montrons ce m\'{e}canisme sur un exemple simple (en annexe, on fournira comme autre exemple la redondance des carr\'{e}s dans le cas du nombre $n=175=5^2.7$).

\noindent La combinatoire \`{a} l'oeuvre, m\^{e}me si elle montre clairement la r\'{e}duction du nombre de petits r\'{e}sidus quadratiques pour les nombres compos\'{e}s, est trop compliqu\'{e}e pour nous permettre de trouver une formule qui donnerait directement ce nombre de petits r\'{e}sidus quadratiques en fonction de $n$.

\noindent $Module \;n=35 \;(R_b(35) \;= \;7 \;\begin{rm}et\end{rm} \;7 \;< \;35/4)$
  
$$
\begin{array}{c|c|c|c|c|c|c|c|c|c|c|c|c|c|c|c|c}
  34 & 33 & 32 & 31 & \cellcolor{cyan}30 & \cellcolor{cyan}29 & 28 & 27 & 26 & \cellcolor{cyan}25 & 24 & 23 & 22 & \cellcolor{cyan}21 & 20 & 19 & 18\\
  \cellcolor{cyan}1 & 2 & 3 & \cellcolor{cyan}4 & 5 & 6 & 7 & 8 & \cellcolor{cyan}9 & 10 & \cellcolor{cyan}11 & 12 & 13 & \cellcolor{cyan}14 & \cellcolor{cyan}15 & \cellcolor{cyan}16 & 17\\
  \hline
  1 & 4 & 9 & 16 & 25 & 1 & 14 & 29 & 11 & 30 & 16 & 4 & 29 & 21 & 15 & 11 & 9\\
  \hline
  \end{array}
$$

\noindent Les redondances de carr\'{e}s pour le module 35 sont :

$$
\begin{array}{rrll}
  6^2 \equiv 1^2 \;(mod \;35) &\begin{rm}car\end{rm} \;(6-1).(6+1)=&5.7 &\begin{rm}et\end{rm} \;35\;|\;35.\\
  11^2 \equiv 4^2 \;(mod \;35) &\begin{rm}car\end{rm} \; (11-4).(11+4)=&7.15=105 &\begin{rm}et\end{rm} \;35\;|\;105.\\
  12^2 \equiv 2^2 \;(mod \;35) &\begin{rm}car\end{rm} \; (12-2).(12+2)=&10.14=105 &\begin{rm}et\end{rm} \;35\;|\;140.\\
  13^2 \equiv 8^2 \;(mod \;35) &\begin{rm}car\end{rm} \; (13-8).(13+8)=&5.21=105 &\begin{rm}et\end{rm} \;35\;|\;105.\\
  16^2 \equiv 9^2 \;(mod \;35) &\begin{rm}car\end{rm} \; (16-9).(16+9)=&7.25=175&\begin{rm}et\end{rm} \;35\;|\;175.\\
  17^2 \equiv 3^2 \;(mod \;35) &\begin{rm}car\end{rm} \; (17-3).(17+3)=&14.20=280 &\begin{rm}et\end{rm} \;35\;|\;280.\\
\end{array}
$$

\noindent Il reste \`{a} \'{e}tudier le cas des produits de nombres premiers simples, de la forme $p.q$ avec $p$ et $q$ premiers (le point (4) \'{e}voqu\'{e} plus haut). L'\'{e}tude sugg\`{e}re que le nombre de petits r\'{e}sidus $R_b(pq)$ est proche de $\displaystyle\frac{R_b(p^2)+R_b(q^2)}{4}$, toujours inf\'{e}rieur \`{a} $\displaystyle\frac{pq}{4}$. Le m\'{e}canisme d'\'{e}limination des carr\'{e}s redondants du fait de l'identit\'{e} remarquable s'applique \`{a} nouveau, faisant diminuer d'autant le nombre de r\'{e}sidus quadratiques.

\noindent Enfin, si l'objectif est de d\'{e}montrer formellement notre assertion en utilisant plut\^{o}t la th\'{e}orie des groupes, peut-\^{e}tre que le fait suivant est \`{a} consid\'{e}rer pour le comptage des petits r\'{e}sidus quadratiques de produits de puissances de nombres premiers diff\'{e}rents : il y a comme un ``triangle'' de produits qui se r\'{e}pondent ; le produit d'un $8k+3$ et d'un $8k'+7$ est un $8k''+5$, celui d'un $8k+3$ et d'un $8k'+5$ est un $8k''+7$ et enfin, le produit d'un $8k+5$ et d'un $8k'+7$ est un $8k''+3$. 

\newpage
\noindent \textbf{Bibliographie}

\noindent [1] Victor-Am\'{e}d\'{e}e Lebesgue, \textit{D\'{e}monstrations de quelques th\'{e}or\`{e}mes relatifs aux r\'{e}sidus et aux non-r\'{e}sidus quadratiques}, Journal de Math\'{e}matiques pures et appliqu\'{e}es (Journal de Liouville), 1842, vol.7, p.137-159.

\vspace{1.7cm}
\noindent \textbf{Annexe 1 : Redondance des carr\'{e}s pour le module $175=5^2.7$}

\noindent Pour le module 175, on \'{e}crit, sous forme de couples, les nombres qui ont m\^{e}me carr\'{e}, on ne pr\'{e}cise pas l'identit\'{e} remarquable $a^2-b^2=(a-b)(a+b)$ qui est telle que les factorisations des nombres $a-b$ et $a+b$ ``contiennent'' tous les facteurs de $175=5^2.7$ :
$$
\begin{array}{l}
  (16,9), (20,15), (23,2), (25,10), (30,5), (32,18), (37,12), (39,11), (40,5), (41,34), (44,19), \\
  (45,10), (46,4), (48,27), (50,15), (51,26), (53,3), (55,15), (57,43), (58,33), (60,10), (62,13), \\
  (64,36), (65,5), (66,59), (67,17), (69,6), (71,29), (72,47), (73,52), (74,24), (75,5), (76,1),\\
  (78,22), (79,54), (80,10), (81,31), (82,68), (83,8), (85,15), (86,61), (87,38).
\end{array}
$$
\noindent De plus, 35 et 70 ont leur carr\'{e} nul et on a pris comme convention de ne pas les compter comme petits r\'{e}sidus quadratiques.

\vspace{1.7cm}

\noindent \textbf{Annexe 2 : petits r\'{e}sidus quadratiques des nombres impairs de 3 \`{a} 51 et leur nombre}

\noindent Entre parenth\`{e}ses est fournie la plus petite racine carr\'{e}e d'un r\'{e}sidu quadratique du module consid\'{e}r\'{e}. 

\begin{tabular}{ll}
\noindent 3 \;$\rightarrow$ 1.&$\rightarrow$ 1\\

\noindent 5 \;$\rightarrow$ 1.&$\rightarrow$ 1\\

\noindent 7 \;$\rightarrow$ 1, 2 (3).&$\rightarrow$ 2\\

\noindent 9 \;$\rightarrow$ 1, 4 (2).&$\rightarrow$ 2\\

\noindent 11 \;$\rightarrow$ 1, 3 (5), 4 (2), 5 (4).&$\rightarrow$ 4\\

\noindent 13 \;$\rightarrow$ 1, 3 (4), 4 (2).&$\rightarrow$ 3\\

\noindent 15 \;$\rightarrow$ 1, 4 (2), 6 (6).&$\rightarrow$ 3\\

\noindent 17 \;$\rightarrow$ 1, 2 (6), 4 (2), 8 (5).&$\rightarrow$ 4\\

\noindent 19 \;$\rightarrow$ 1, 4 (2), 5 (9), 6 (5), 7 (8), 9 (3).&$\rightarrow$ 6\\

\noindent 21 \;$\rightarrow$ 1, 4 (2), 7 (7), 9 (3).&$\rightarrow$ 4\\

\noindent 23 \;$\rightarrow$ 1, 2 (5), 3 (7), 4 (2), 6 (11), 8 (10), 9 (3).&$\rightarrow$ 7\\

\noindent 25 \;$\rightarrow$ 1, 4 (2), 6 (9), 9 (3), 11 (6).&$\rightarrow$ 5\\

\noindent 27 \;$\rightarrow$ 1, 4 (2), 7 (13), 9 (3), 10 (8), 13 (11).&$\rightarrow$ 6\\

\noindent 29 \;$\rightarrow$ 1, 4 (2), 5 (11), 6 (8), 7 (6), 9 (3), 13 (10).&$\rightarrow$ 7\\

\noindent 31 \;$\rightarrow$ 1, 2 (8), 4 (2), 5 (6), 7 (10), 8 (15), 9 (3), 10 (14), 14 (13).&$\rightarrow$ 9\\

\noindent 33 \;$\rightarrow$ 1, 3 (6), 4 (13), 9 (3), 12 (12), 15 (9), 16 (4).&$\rightarrow$ 7\\

\noindent 35 \;$\rightarrow$ 1, 4 (2), 9 (3), 11 (9), 14 (7), 15 (15), 16 (4).&$\rightarrow$ 7\\

\noindent 37 \;$\rightarrow$ 1, 3 (15), 4 (2), 7 (9), 9 (3), 10 (11), 11 (14), 12 (7), 16 (4).&$\rightarrow$ 9\\

\noindent 39 \;$\rightarrow$ 1, 3 (9), 4 (2), 9 (3), 10 (7), 12 (18), 13 (13), 16 (4).&$\rightarrow$ 8\\

\noindent 41 \;$\rightarrow$ 1, 2 (17), 4 (2), 5 (13), 8 (7), 9 (3), 10 (16), 16 (4), 18 (10), 20 (15).&$\rightarrow$ 10\\

\noindent 43 \;$\rightarrow$ 1, 4 (2), 6 (7), 9 (3), 10 (15), 11 (21), 13 (20), 14 (10), 15 (12), 16 (4), 17 (19), 21 (8).&$\rightarrow$ 12\\

\noindent 45 \;$\rightarrow$ 1, 4 (2), 9 (3), 10 (10), 16 (4), 19 (8).&$\rightarrow$ 6\\

\noindent 47 \;$\rightarrow$ 1, 2 (7), 3 (12), 4 (2), 6 (10), 7 (17), 8 (14), 9 (3), 12 (23), 14 (22), 16 (4), 17 (8), 18 (21), 21 (16).&$\rightarrow$ 14\\

\noindent 49 \;$\rightarrow$ 1, 2 (10), 4 (2), 8 (20), 9 (3), 11 (16), 15 (8), 16 (4), 18 (19), 22 (13), 23 (11).&$\rightarrow$ 11\\

\noindent 51 \;$\rightarrow$ 1, 4 (2), 9 (3), 13 (8), 15 (24), 16 (4), 18 (18), 19 (11), 21 (15), 25 (5).&$\rightarrow$ 10\\
\end{tabular}

\end{document}